\documentclass[11pt]{amsart}

\usepackage[a4paper,margin=1in]{geometry}
\usepackage{amsmath,amssymb,amsthm,mathtools}
\usepackage{microtype}

\usepackage[english]{babel}
\usepackage{psfrag, url}
\usepackage[latin1]{inputenc}
\usepackage[dvips]{graphicx}
\usepackage{textcomp}
\usepackage{flafter}

 \usepackage{esint}
 
\usepackage{epsfig}
\usepackage{bbm}
\usepackage{wrapfig}
 
\usepackage{eurosym}
\usepackage[hidelinks]{hyperref}  

 \theoremstyle{plain}
\newtheorem{thm}{Theorem} 
\newtheorem{cor}[thm]{Corollary}

\newtheorem{conj}[thm]{Conjecture}

\theoremstyle{definition}

\theoremstyle{remark}
\newtheorem{rem}[thm]{Remark}

\newcommand{\R}{\mathbb{R}}

\newcommand{\cH}{\mathcal{H}}

\newcommand{\Sing}{\mathrm{Sing}}
\newcommand{\Reg}{\mathrm{Reg}}

\newcommand{\res}{\mathop{\hbox{\vrule height 7pt width .5pt depth 0pt
			\vrule height .5pt width 6pt depth 0pt}}\nolimits}

 \newcommand{\eps}{\varepsilon}
  \newcommand{\cJ}{\mathcal{J}}

\title[The classical Bernoulli free boundary problem]{\small Overview and recent progress on\\ the classical Bernoulli free boundary problem}
\author{Xavier Fern\'andez-Real}
\thanks{The author was supported by the Swiss National Science Foundation (SNF grant PZ00P2\_208930), by the Swiss State Secretariat for Education, Research and Innovation (SERI) under contract number MB22.00034, and by the AEI project PID2024-156429NB-I00 (Spain).}

\date{}

\begin{document}

\begin{abstract}
    This survey gives an overview of the Bernoulli or one-phase free boundary problem, with emphasis on recent progress. 
\end{abstract}

\maketitle

\tableofcontents

\section{Introduction}

Free boundary problems are variational or PDE models in which part of the boundary of the relevant domain is \emph{unknown} and must be determined together with the solution.
In contrast to boundary value problems posed in a fixed region, here the geometry of the interface is an additional unknown, and the equations typically couple bulk behavior to boundary conditions \emph{on the a priori unknown set}.
This coupling is a natural way to encode ``phase'' separation: one phase is described by a PDE in an open set, while the complement is constrained by an obstacle/contact condition, and the interface between them is the \emph{free boundary}.

In recent years, one of the main drivers of development in the theory of nonlinear elliptic
PDE has been its connection with physical models involving phase transitions. This line of study
has been of particular mathematical interest because of the interdisciplinarity of its
techniques, and because, while it was originally nurtured by theories developed in other areas,
the flow of information has by now become genuinely bi-directional. Free boundary problems are a paradigmatic example of such a connection: at the intersection of Geometric Analysis and Partial Differential Equations, their study fully benefits from techniques coming from both areas. Furthermore, especially in the last few years, we have seen how some of the insights obtained in the free boundary setting can be translated back to the aforementioned contexts.
For instance, recent proofs of one-dimensional symmetry for bounded stable solutions of the classical Allen--Cahn equation in $\R^3$ build on ideas developed in the free-boundary setting \cite{chan2026stableallen,liu2026stabledegiorgi}. 

In this expository note we study one of the most classical free boundary problems, and arguably the main example of a nonconvex free boundary problem: the \emph{one-phase Bernoulli problem}, also called the \emph{Alt--Caffarelli problem}\footnote{In this work, we use the names one-phase problem, Bernoulli problem, and Alt--Caffarelli problem interchangeably}. In its classical form, it corresponds to studying the Serrin-type problem of searching for a nonnegative function $u$ and an open set $\{u > 0\} \subset \R^d$ such that, in a given domain $D$ (with possibly $D = \R^d$) we have
\begin{equation}
\label{eq:onephase0}
\left\{
\begin{array}{rcll}
\Delta u &=& 0&\quad\text{in}\quad \{u > 0\}\cap D,\\
u & \ge & 0 & \quad\text{in}\quad D, \\
|\nabla u| & = & 1 & \quad\text{on}\quad \partial \{ u > 0\}\cap D,
\end{array}
\right.
	\end{equation}
where the boundary condition $|\nabla u|=1$ is understood in an appropriate sense (classical/viscosity, depending on the regularity available).
The last condition is \emph{overdetermined} from the viewpoint of a fixed domain: it prescribes the magnitude of the normal derivative on the boundary (exactly where $u= 0$ as well), and it is precisely this extra condition that selects the geometry of $\partial\{u>0\}$. Together with the solution, the unknown domain that must be determined is, in this case, the positivity set of the nonnegative function $u$, $\{u > 0\}$.
 
 Alt and Caffarelli investigated the behaviour of \eqref{eq:onephase0} by observing that it corresponds to the criticality condition of what is currently known as the \emph{Alt--Caffarelli functional}: 
 \begin{equation}
 \label{eq:onephase}
 \cJ (u) = \int_{D}\left( |\nabla u|^2 + \mathbbm{1}_{\{u > 0\}}\, \right)\, dx = \int_{D}|\nabla u|^2\, dx + |\{u > 0\}\cap D| ,\quad u:D \to [0, +\infty).
 \end{equation}

The one-phase Alt--Caffarelli free boundary problem is a prototypical model for a variety
of interface phenomena. It arises in idealized fluid-mechanical regimes (notably jet flows
and cavitational flow), in combustion-related models such as flame propagation, and electrostatics
\cite{BL82,CV95,PY07,kriventsov2025min,sirakov2002overdetermined}. It also appears naturally in optimal design and shape optimization,
where free boundaries encode the geometry selected by competing bulk and boundary
effects \cite{AC81,KL18,KL19}. From a mathematical point of view, the variational formulation
was introduced by Alt and Caffarelli in \cite{AC81}, and the subsequent works of Caffarelli and
collaborators developed a far-reaching regularity and structure theory for minimizers and
related solution classes \cite{Caf87,Caf89,Caf88,CJK04}. These contributions are among the
results recognized by Caffarelli's Wolf Prize (2012) and Abel Prize (2023). The subject
remains very active, with many recent developments in analysis and geometry of free
boundaries \cite{DJ09,JS15,EE19,ESV20,DSV21,EFeY23,FY23,ESV24, openai2026critical}; see also the
monographs \cite{CS05,Vel23}.

Historically, the problem was motivated by questions from fluid mechanics \cite{AC81}.
In that setting, the Bernoulli condition models interfaces in idealized flow regimes; even in
simplified configurations, the geometry of the interface can be remarkably subtle.
A classical manifestation of this flexibility is the existence of nontrivial critical points with
intricate geometry related to vortex arrays and ``hollow vortex'' configurations \cite{BSS76}.
The one-phase framework also appears in other evolution and combustion models through
quasi-static limits or associated free boundaries; for instance, free boundary problems arise
naturally in flame propagation models \cite{CV95}.

A different but closely related source is the theory of \emph{overdetermined elliptic problems}.
Prescribing both Dirichlet and Neumann-type boundary data typically forces rigidity
(Serrin-type phenomena, where under appropriate assumptions the domain must be a ball),
but relaxing the hypotheses (e.g.\ allowing unbounded or periodic domains, or adopting weaker
notions of solution) leads to rich families of ``exceptional'' domains and reconnects directly
with one-phase free boundaries; see, for example, \cite{HHP11,Tra14, liu2021smooth} and the references therein.

Finally, the one-phase or Bernoulli problem is tightly linked to semilinear elliptic equations and
phase-transition models through blow-down and singular perturbation limits. For instance,
in the Alt--Phillips functional \cite{phillips1983minimization, AP86}, the volume term in \eqref{eq:onephase} is
replaced by a power of $u$, producing a family of free boundary problems whose
regularity and stability theory parallels the Alt--Caffarelli setting while introducing new
analytic features. A complementary connection, closer to combustion models, comes from
bump-type semilinear equations. If
\[
\mathcal E_\eps(v,\Omega) = \int_\Omega \bigl(|\nabla v|^2+\Phi(v/\eps)\bigr)\,dx,
\qquad \Phi'=\beta,\quad \Phi(0)=0,\quad \beta\ge0,\quad \int_0^\infty \beta=1,
\]
then critical points solve
\[
\Delta v_\eps=\frac1{2\eps}\beta(v_\eps/\eps),
\]
and, as $\eps\downarrow0$, the potential term formally converges to
$\chi_{\{v>0\}}$ (and the right-hand side in the PDE formally converges to a measure concentrated on the free boundary). Thus the thin reaction region $\{v_\eps\sim\eps\}$ becomes, in the
limit, the Bernoulli free boundary. Conversely, blow-downs of global Lipschitz solutions
of equations of the form $\Delta u=f(u)$, with $f$ of bump type, lead to global
one-phase solutions \cite{fr2019dcds,audrito2022interface}. In this sense the one-phase problem plays, for
these nonbalanced semilinear equations, a role analogous to the role played by minimal
surfaces in the Allen--Cahn theory. This viewpoint connects free-boundary regularity and
rigidity to De Giorgi/Bernstein-type questions for semilinear equations, and has recently
been used to obtain one-dimensional symmetry results for semilinear minimizers and
monotone solutions \cite{fr2019dcds,audrito2022interface,EFeY23}.

\bigskip

 Minimizers (and more generally, critical points) of \eqref{eq:onephase} (with prescribed boundary data) satisfy $\Delta u=0$ in 
 \[
 \Omega_u := D \cap \{u>0\}
 \]
 and a Bernoulli-type condition $|\nabla u| = 1$ on 
 \[
 \Gamma_u:=D\cap \partial\{u>0\};
 \] 
 informally, variations that change the positivity set balance Dirichlet energy against the volume term and produce a Bernoulli condition with constant gradient magnitude along $\Gamma_u$.
This variational structure places the problem at the intersection of elliptic PDE, geometric measure theory, and calculus of variations; it also makes blow-up analysis and monotonicity tools available.

\begin{samepage}
For the sake of presentation, in this note, we will always consider \emph{solutions}\footnote{The notion of solution is, at the moment, left intentionally vague. A solution could refer to a minimizer of \eqref{eq:onephase}, to a viscosity solution, to a stable solution with or without density bounds, to a stationary solution, etc. For the reader's convenience, the simplest way to think about it and not have to deal with regularity issues at this point is to consider \emph{classical solutions} to \eqref{eq:onephase0}: those where the free boundary $\Gamma_u$ is smooth. We are then interested in estimates that are independent of its regularity or smoothness. } to \eqref{eq:onephase0} that are locally Lipschitz continuous\footnote{In some cases, like for minimizers, classical solutions, or viscosity solutions, being Lipschitz continuous is not a hypothesis but rather a property. In dimension two, local Lipschitz regularity also follows from inner stationarity for nonnegative scalar $H^1$ functions,  see \cite[Corollary~7]{dephilippis2026optimal}. } in the domain and belong to $H^1_{\rm loc}(D)$, 
\[
u\in C^{0, 1}_{\rm loc}(D)\cap H^1_{\rm loc}(D). 
\]
\end{samepage}

We also take this opportunity to observe that the problem is scale invariant under dilations: if $u$ satisfies \eqref{eq:onephase0} in a domain $D$ and $r > 0$, then the rescalings 
\begin{equation}
\label{eq:urrescaling}
u_r(x) := \frac{u(rx)}{r},
\end{equation}
satisfy the same equation, \eqref{eq:onephase0}, in 
\[
D_r := \frac{1}{r}D. 
\]

\section{Minimizers}

The most studied class of solutions to the one-phase problem is given by minimizers of \eqref{eq:onephase}. Namely, given $D\subset \R^d$ a bounded domain, we consider elements of the set of minimizers: 
\[
\left\{ u\in C^{0, 1}_{\rm loc}(D)\cap H^1(D) : \cJ(u) \le \cJ (v) \ \  \forall v\in H^1(D)\ \ \text{with}\ \ v - u \in H^1_0(D) \right\}.
\]

On the one hand, it is natural to start by understanding absolute minimizers of the functional: they are the most energetically efficient configuration, and they can always be constructed or found regardless of ($H^1$) boundary data. On the other hand, minimizers enjoy structural and rigidity properties that make them much more tractable than other non-minimizing critical points or simply solutions to \eqref{eq:onephase0}. Thanks to that, from a regularity point of view, minimizers have a particularly favourable behaviour. 

For instance, if $u$ is a minimizer of $\cJ$ in $B_1\subset \R^d$, and $0\in \Gamma_u$, then $u$ satisfies the growth condition
\begin{equation*}
 \fint_{\partial B_r} u \, d\cH^{d-1} \le C_\circ r\qquad\text{for all}\quad 0  < r < 1/2. 
\end{equation*}
for some $0<C_\circ <+\infty$, and is \emph{nondegenerate}:
\begin{equation*}
 \fint_{\partial B_r} u \, d\cH^{d-1} \ge \delta_\circ r\qquad\text{for all}\quad 0  < r < 1/2, 
\end{equation*}
for some $\delta_\circ > 0$. Moreover, a minimizer $u$ satisfies the density estimates: 
\begin{equation}
    \label{eq:density_estimates}
    0 < \delta_\circ \le \frac{|\{u > 0\}\cap B_r|}{|B_r|}\le 1-\delta_\circ\qquad\text{for all}\quad 0  < r < 1/2. 
\end{equation}
The most comprehensive result on the regularity of free boundaries for minimizers of the Alt--Caffarelli functional is the following (see \cite{CJK04, JS15, Wei99} for results up to $d = 4$, and the recent OpenAI preprint \cite{openai2026critical} for the remaining dimensions $d = 5, 6$): 
\begin{thm}
    Let $d \ge 2$, and let $D\subset \R^d$ be a bounded domain. Let $u\in H^1(D)$ be a minimizer of \eqref{eq:onephase} such that $u - g \in H^1_0(D)$ for some $0\le g\in H^1({D})$. Then, we can split 
    \[
    \Gamma_u = \Reg(u) \cup \Sing(u), 
    \]
    where $\Reg(u)$ is relatively open, locally a $C^\infty$ (analytic) $(d-1)$-dimensional manifold, and: 
    \begin{itemize}
        \item For $d \le 6$, $\Sing(u) = \varnothing$. 
        \item For $d \ge 7$, ${\rm dim}_{\cH} \left( \Sing(u) \right) \le d-7$, where ${\rm dim}_\cH$ denotes the Hausdorff dimension of a set. 
    \end{itemize}
\end{thm}

\begin{figure}[t]
    \includegraphics[scale = 1.1]{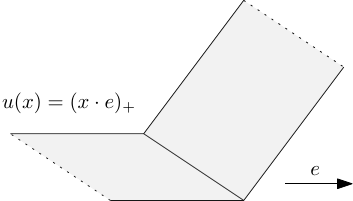}
    \centering
    \caption{\small The canonical solution to the one-phase problem is the half-space solution, $u(x) = (x\cdot e)_+$ for some $e\in \mathbb{S}^{d-1}$.}
    \label{fig:HalfSpace}
\end{figure}

The result is in a way reminiscent of that in minimizing minimal surfaces: for low dimensions, the manifolds being studied (free boundaries or hypersurfaces) happen to be smooth, whereas in large dimensions there may be singularities. Moreover, the canonical example of global solutions in both cases is the hyperplane, which in the case of the Bernoulli problem corresponds to the free boundary of the solution\footnote{The hyperplane also appears as the free boundary of $|x\cdot e|$, which in a way would correspond to a multiplicity 2 hyperplane in the case of minimal surfaces. In this case, however, the solution is not minimizing.} 
\begin{equation}\label{eq:half-space}
u(x) = (x\cdot e)_+\quad\text{for some}\quad e\in \mathbb{S}^{d-1}. 
\end{equation}
This is known as a \emph{half-space solution} (see Figure~\ref{fig:HalfSpace}). Furthermore, the previous result can be made quantitative via a compactness argument (cf. \cite{CFFS25}):
\begin{cor} 
\label{cor:minimizers_regularity}
    Let $2\le d \le 6$, and let  $u\in H^1(B_1)$ be a minimizer of \eqref{eq:onephase} with $0\in \partial \{u > 0\}$. Then, 
    \[
    |D^2 u |\le C \quad\text{in}\quad \{u > 0\}\cap B_{1/2}, 
    \]
    for some universal constant $C$. 
    In particular, the second fundamental form of the free boundary $|{\rm II}_{\Gamma_u}|$ is bounded by a universal constant in $B_{1/2}$, and the contact set $\{u = 0\}$ can be locally written as the subgraph of a smooth function in universal balls around each free boundary point. 
\end{cor}

The analogies with minimal surfaces are far from ending here, though. The strategy to prove the previous theorem follows very closely the one employed in the case of minimal surfaces as well, even if the particularities of each step differ: 
\begin{enumerate}
    \item Up to a translation, we may assume $0\in \Gamma_u$. The natural homothetic rescaling $u_r$ from \eqref{eq:urrescaling} keeps the problem invariant, with $0 \in \Gamma_{u_r}$. 
    \item The problem enjoys a monotonicity formula, known as Weiss' monotonicity formula. This means that rescalings have a monotone energy (with respect to the scale), which can be used to control how close a solution is to its 1-homogeneous extension. 

    More precisely, if $u$ is a minimizer of the Alt-Caffarelli functional in $B_1$ (or more generally, a stationary solution, see \eqref{eq:stationary_def}) with $0\in \Gamma_u$, then the function
    \[
    r\mapsto W(u, r) := r^{-d}\int_{B_r}\left( |\nabla u|^2 + \mathbbm{1}_{\{u > 0\}}\right)\,  dx-r^{-d-1} \int_{\partial B_r} u^2\, d\mathcal{H}^{d-1}
    \]
    is nondecreasing, and it is constant if and only if $u$ is 1-homogeneous. 
    \item Coupled with the fact that we can control the Lipschitz regularity of solutions, one may perform a blow-up argument around 0 letting $r\downarrow 0$ for $u_r$ (along a subsequence) to obtain a global one-homogeneous solution (minimizer), $u_\infty$. 
    \item In parallel, solutions to the one-phase problem enjoy an ``improvement-of-flatness'' type result: any solution close enough (in some appropriate sense) to the canonical solution \eqref{eq:half-space} necessarily has a smooth free boundary in a smaller domain. Namely, if the free boundary is sufficiently flat, then it is smooth. 
    \item It only remains to classify possible blow-ups $u_\infty$. The only 1-homogeneous  smooth free boundary is the hyperplane, corresponding to \eqref{eq:half-space}. If one can say that \emph{all} 1-homogeneous minimizers are of the form \eqref{eq:half-space}, then we necessarily have that solutions have smooth free boundaries. Up to dimension 6, we know this is the case, with counter-examples in higher dimensions.
\end{enumerate}

\subsection{Classification of 1-homogeneous minimizers}\label{ssec:classification_homogeneous}

As explained in the previous discussion, the regularity theory for the free boundary of minimizers is reduced (via blow-up) to the classification of global $1$-homogeneous minimizers. In particular, singular free-boundary points correspond to nonflat $1$-homogeneous blow-ups, so understanding whether such cones exist is the key structural question.

A useful point of view is that the classical classification results through dimension six are proved at the level of \emph{stable critical points} with density estimates (or smooth outside of the origin), which is a larger class than minimizers (recall \eqref{eq:density_estimates}).

Let $u$ be a global $1$-homogeneous solution of the one-phase problem in $\R^d$, and assume that $\Gamma_u$ is smooth outside of the origin\footnote{This hypothesis is automatic if one argues in the smallest dimension where a nonflat singular cone exists, by a Federer-type dimension reduction argument.}, in the sense of Remark~\ref{rem:smoothoutside} below. We say that $u$ is a \emph{critical point} (with respect to inner variations, cf. \eqref{eq:stationary_def} below for general non-homogeneous solutions) if
\[
\left.\frac{d}{dt}\right|_{t=0}\cJ\big(u(\,\cdot+tX(\cdot)\,)\big)=0
\]
for every smooth domain $D\Subset \R^d\setminus\{0\}$ and every smooth vector field $X$ compactly supported in $D$ (where $\cJ$ is computed in $D$). We say that $u$ is \emph{stable} if, in addition,
\[
\left.\frac{d^2}{dt^2}\right|_{t=0}\cJ\big(u(\,\cdot+tX(\cdot)\,)\big)\ge 0
\]
for every such $D$ and $X$. Of course, every minimizing cone is a stable critical point.

For smooth cones (that is, 1-homogeneous solutions that are smooth outside of the origin), the second variation can be rewritten as the (Caffarelli--Jerison--Kenig, \cite{CJK04}) stability inequality
\[
\int_{\{u>0\}} |\nabla \xi|^2\,dx \;\ge\; \int_{\Gamma_u} H\,\xi^2\,d\cH^{d-1}
\qquad\text{for all }\xi\in C_c^1(\R^d\setminus\{0\}),
\]
where $\nu=\nabla u$ is the inward unit normal to $\{u>0\}$ and $H=-u_{\nu\nu}$ denotes the mean curvature of $\Gamma_u$ with this convention; for these cones, $H\ge0$.

The first positive low-dimensional classification result for minimizers is then the following, which  we more generally state for stable solutions with positive density of the zero set\footnote{In the statement of the result, the density estimates are alternatively enforced by assuming that the cone is smooth outside of the origin.}.

\begin{thm}[\cite{CJK04, JS15, openai2026critical}]
Let $u$ be a nonzero $1$-homogeneous stable critical point of the Alt--Caffarelli functional in $\R^d$, and assume that $\Gamma_u$ is smooth outside of the origin. If $d\le 6$, then $u$ is one-dimensional, i.e.
\[
u(x)=(x\cdot e)_+ \qquad\text{for some }e\in \mathbb{S}^{d-1}.
\]
\end{thm}

\begin{rem}\label{rem:smoothoutside}
    Being smooth outside of the origin means the following: for each free boundary point $x_\circ$ different from the origin, the contact set (or alternatively, the positivity set) can locally be written as the subgraph of a smooth function around $x_\circ$. In particular, smooth free boundaries cannot have vanishing density of the contact set, and solutions of the form $|x\cdot e|$ are \emph{not} smooth outside of the origin. 
\end{rem}

Here the case $d=3$ is due to Caffarelli--Jerison--Kenig, \cite{CJK04}, while Jerison--Savin proved the case $d=4$, \cite{JS15}. The cases $d=5,6$ follow from the very recent preprint by OpenAI, \cite{openai2026critical}, which refines the strategy of Jerison--Savin by coupling a scalar test with a vector field, built from the normalized Hessian, to obtain compatible interior and boundary inequalities. In particular, since minimizers are stable, this yields the classification of $1$-homogeneous minimizing cones in dimensions $d\le 6$.

\begin{figure}[t]
    \includegraphics[scale = 0.65]{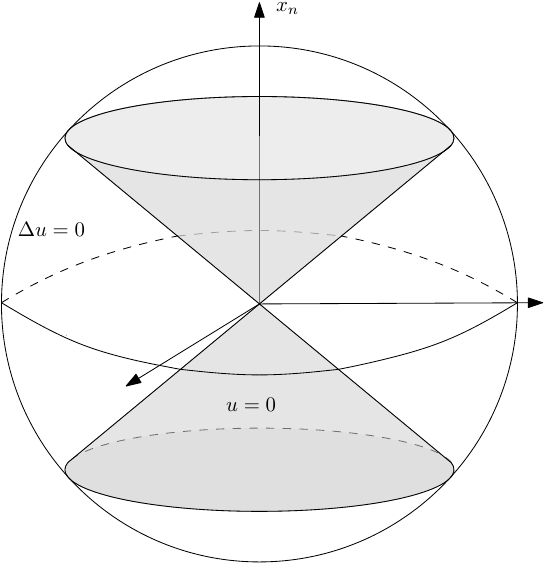}
    \centering
    \caption{\small The De Silva--Jerison cone is an axially symmetric solution to the Bernoulli problem, which can be constructed in all dimensions $d \ge 3$ for a suitable aperture and is stationary in each such dimension (see \eqref{eq:stationary_def}). In dimension $d=7$, De Silva and Jerison proved that this solution is minimizing \cite{DJ09}.}
    \label{fig:DeJ}
\end{figure}

The Jerison--Savin argument is to prove that every nonflat homogeneous cone is \emph{unstable}. They do so by constructing a positive homogeneous function \(\bar v\) which is a \emph{strict subsolution} of the linearized free-boundary problem outside the origin. The construction is based on taking an appropriate power of a quantity \(w=w(D^2u)\): one proves a Simons-type \emph{interior} differential inequality for \(w\) (equivalently, for \(\log w\)) and a complementary \emph{boundary} inequality on \(\Gamma_u\), obtained from the identities satisfied by \(D^2u\) at free-boundary points. In dimension \(3\), one can take \(w=|D^2u|\). In dimension \(4\), the choice \(w=|D^2u|\) does not yield a strong enough boundary estimate, and the key new idea is to replace it by a weighted expression in the eigenvalues of \(D^2u\),
\[
w^2=\sum_{\lambda_k>0}\lambda_k^2+4\sum_{\lambda_k<0}\lambda_k^2,
\]
which still satisfies the same interior inequality (by convexity) but improves the boundary inequality enough to produce the required strict subsolution. This gives the contradiction to stability, and hence every stable \(1\)-homogeneous cone in \(d\le 4\) must be the half-space solution.

The picture changes in dimension $7$.

\begin{thm}[\cite{DJ09}]\label{thm:desilvajerison}
There exists a nonflat $1$-homogeneous minimizer $U$ in $\R^7$ whose free boundary is a cone and is smooth outside of the origin. In particular, $0$ is an isolated singular free boundary point.
\end{thm}

Their example is the axially symmetric conical \emph{critical point} \(U\) (the free-boundary analogue of the Simons cone in minimal-surface theory, and often called the \emph{De Silva--Jerison cone}, see Figure~\ref{fig:DeJ}). The proof that \(U\) is minimizing is a comparison/trapping argument: one fixes \(U\) as boundary datum on a ball, lets \(u\) be an energy minimizer with those boundary values, and then constructs a strict subsolution \(V\le U\) and a strict supersolution \(W\ge U\). By dilating these barriers, one obtains families \(V_\rho\) and \(W_\rho\) that approach \(U\) from below and above and can be compared with \(u\); letting the dilation parameter tend to infinity forces
\[
V_\rho \le u \le W_\rho \quad\Longrightarrow\quad u=U.
\]
A key point is that the sub- and supersolution constructions are genuinely asymmetric (the supersolution is much more delicate), and a substantial part of the verification reduces to explicit estimates for Legendre-type functions, which are checked with computer assistance.

It is worth stressing that, in the same axial symmetry class, Caffarelli--Jerison--Kenig had already derived a criterion that rules out the nonflat axially symmetric cone as a minimizer in dimensions $d\le 6$. Thus the De Silva--Jerison cone appears exactly in the first dimension not excluded by that argument. We also refer to Subsection~\ref{sec:more_singular_cones} below for more recent examples of singular cones. 

Thus, the critical dimension $d_\ast$ (the first dimension with singular minimizing cones) satisfies
\[
d_\ast=7.
\]

\subsubsection{Generic regularity}
There is also a \emph{generic} version of this picture, which is now known.
For admissible monotone one-parameter families of boundary data, \cite{FY23} and
\cite{fernandez2025linearized} extend the range of generic regularity
by one and two dimensions, respectively, for generic boundary data
(in an appropriate sense). In particular, since $d_\ast=7$, the free boundaries of minimizers are smooth for almost every parameter in such a family in dimensions $d\le8$.

\subsubsection{Structure near singular points}

The classification of cones determines which singularities may occur, but there is a
finer local question: once a singular cone appears as a blow-up, does the solution
actually approach it in a unique and quantitative way? Weiss' monotonicity formula
implies that every blow-up of a minimizer at a free-boundary point is a global
$1$-homogeneous minimizer, but by itself it does not rule out different tangent cones
along different sequences of scales \cite{Wei99}.

At the level of the singular set, quantitative stratification gives a robust answer:
the singular set of Alt--Caffarelli minimizers is rectifiable and satisfies
quantitative measure bounds \cite{EE19}. Thus, singular points organize into lower-dimensional rectifiable strata (analogously to dimension-reduction and quantitative-stratification theory for minimal surfaces and harmonic maps).

At individual singular points, the sharpest description is known when the tangent cone
has an isolated singularity. If one blow-up $U$ of a minimizer $u$ at
a point in ${\rm Sing}(u)$ is such a minimizing cone, then Engelstein--Spolaor--Velichkov prove
that $U$ is the unique blow-up and that 
\[
\partial\{u>0\}
\quad\text{is a } C^{1,\log}\text{ graph over }\partial\{U>0\}
\]
at small scales. If the cone is integrable through rotations, the convergence improves
to $C^{1,\alpha}$; in particular this applies to the De Silva--Jerison cones
\cite{ESV20}. (See also \cite{engelstein2025asymptotic} for a more recent, similar result without requiring minimality; and \cite{edelen2024symmetric} for a symmetric $(\log)$-epiperimetric perspective on decay and growth of the Weiss energy.) Thus, near the known isolated minimizing singularities, the free boundary
has a genuine conical expansion, not merely subsequential conical tangents.

\subsubsection{Further examples of singular minimizing cones}\label{sec:more_singular_cones}
Finally, the recent work of Firester--Tsiamis--Wang, \cite{firester2026stability, firester2026area} shows that the landscape of explicit singular minimizing cones is much richer than previously known. In a companion paper (via a capillary deformation/barrier construction refining the De Silva--Jerison strategy), they construct many new explicit singular minimizing one-phase cones in high dimensions: for instance, they obtain a second minimizing cone in dimension $7$ (besides the De Silva--Jerison example), and more generally prove global minimality for a large family of bi-orthogonally symmetric cones $U_{d,k}$ in dimensions $d\ge 7$.

Further examples have been obtained by Carducci--Velichkov, who prove that the Cartan one-phase cone in $\R^{26}$ is minimizing, providing an example outside the axially symmetric and Lawson-type families; the corresponding Cartan cones in dimensions $5,8,14$ are unstable \cite{carducci2026minimizing}.

\section{Stable solutions}

\label{sec:stable}
Beyond minimizers, an important class of solutions is given by \emph{stable} critical points of the Alt--Caffarelli functional (without imposing density constraints). This is the natural free-boundary analogue of stable minimal surfaces (or stable critical points for semilinear elliptic energies), and it is broad enough to include many geometrically relevant examples, while still retaining useful rigidity and compactness properties. It is also the class of physically observable configurations: in a natural setting, there is no reason to expect a solution to be a global minimizer (and in practice, this is very hard to verify), whereas one does expect an observed configuration to be at least locally stable, in the sense that its energy does not immediately decay under small perturbations. This makes stability a more flexible and robust notion than minimality, and it also raises the question of (quantitative) regularity of free boundaries in this setting. A systematic treatment of this class has been developed recently in \cite{kamburov2022nondegeneracy, CFFS25, fernandezreal2026stability}.

Let $D\subset \R^d$ be a domain. If $u\in H^1_{\rm loc}(D)$ and $X\in C_c^\infty(D;\R^d)$, we consider inner variations (as in the case for cones in Subsection~\ref{ssec:classification_homogeneous}) and we say that $u$ is a \emph{stationary solution} (or stationary critical point) in $D$ if
\begin{equation}
\label{eq:stationary_def}
\left.\frac{d}{dt}\right|_{t=0}\cJ\big(u(\cdot + t X(\cdot) )\big)=0
\qquad \text{for all }X\in C_c^\infty(D;\R^d);
\end{equation}
and we say that $u$ is \emph{stable} in $D$ if it is stationary and
\[
\left.\frac{d^2}{dt^2}\right|_{t=0}\cJ\big(u(\cdot + t X(\cdot))\big)\ge 0
\qquad \text{for all }X\in C_c^\infty(D;\R^d).
\]
This is the inner-variation notion of stability. In parallel, one keeps track of the more geometric class of \emph{classical solutions}, namely nonnegative functions $u$ such that $\{u>0\}$ has smooth free boundary in $D$ and satisfies \eqref{eq:onephase0}. A \emph{classical stable solution} is one that is both classical and stable in the above sense. This is the setup used in \cite{kamburov2022nondegeneracy, CFFS25, fernandezreal2026stability}, for example.

\begin{figure}[t]
    \includegraphics[scale = 1.0]{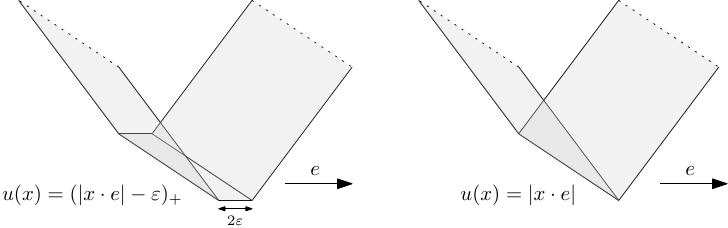}
    \centering
    \caption{\small The vee solution can be obtained as a limit of classical stable solutions.}
    \label{fig:veeSol}
\end{figure}

It is important not to work with arbitrary stationary/stable inner-variation critical points only. There are spurious examples (for instance $u(x_1,x_2)=|x_1x_2|$ in $\R^2$) that are stationary and stable under domain variations, but do not behave like genuine free-boundary solutions and should not be part of the regularity theory. For this reason, a robust class is given by \emph{local limits of classical stable solutions}: many regularity and rigidity statements extend to this class, and this is essentially the largest class where one expects them to remain true. In this exposition we will call \emph{stable solutions} those that are either \emph{classical critical points} with nonnegative second variation, or limits of solutions in this class.

A useful point to keep in mind is that stability is strictly weaker than minimality. For example, the canonical ``vee'' solution
\[
u(x)=|x\cdot e|
\]
is a global stable solution whose free boundary is a hyperplane, and it is a limit of classical solutions ($(|x\cdot e|-\eps)_+$ as $\eps\downarrow 0$, see Figure~\ref{fig:veeSol}), but it is not minimizing for the one-phase functional (it is a free-boundary analogue of a multiplicity-two hyperplane). More generally, one can glue minimizing configurations along a common zero/contact region to produce stable, non-minimizing solutions. This phenomenon is one of the reasons why the stable class is richer, and harder, than the minimizing class.

For classical stable solutions, the second variation can be rewritten in a geometric form as in the version of Caffarelli-Jerison-Kenig above, so that it is equivalently written as the nonnegativity of the following second variation quadratic form on compactly supported variations
\begin{equation}
\label{eq:stab1}
\mathcal{Q}(\xi,\xi) := \int_{D\cap \{u > 0\}} |\nabla \xi|^2\,dx -  \int_{D\cap \Gamma_u} H\,\xi^2\, d\mathcal{H}^{d-1}
\ge 0 \qquad \forall  \xi \in C_c^\infty\big(D \cap \overline{\{u > 0\}}\big),
\end{equation} 
where $\nu=\nabla u$ is the inward unit normal to $\{u>0\}$ and $H=-u_{\nu\nu}$ denotes the mean curvature of the free boundary $\Gamma_u$ with this convention. This in turn yields a Sternberg--Zumbrun type bulk inequality:
\begin{equation}
    \label{eq:stab2}
\int_{D\cap\{u>0\}}
\bigl(|D^2u|^2-|\nabla |\nabla u||^2\bigr)\eta^2 \le \int_{D\cap\{u>0\}} |\nabla u|^2 |\nabla \eta|^2 \qquad\text{for all }\eta\in C_c^1(D),
\end{equation}
which is particularly useful in certain applications. A key consequence is the scale-invariant Hessian control
\[
\int_{B_R\cap\{u>0\}} |D^2u|^2 \le C R^{d-2}
\]
for global classical stable solutions (with $C$ depending only on the dimension). This estimate plays a central role in compactness and blow-down arguments. 

As in the minimizing theory, blow-ups and blow-downs naturally lead to homogeneous solutions (which do not vanish because stable solutions enjoy nondegeneracy as well, \cite{kamburov2022nondegeneracy}), so the classification of stable $1$-homogeneous cones is a central issue. Critically, we observe that contrary to the case of minimizers, the vee solution is always a stable solution in all dimensions. This makes the path towards a regularity theory described above for minimizers (e.g., Corollary~\ref{cor:minimizers_regularity}) much more complicated: the improvement of flatness cannot always be applied, since one needs to consider solutions that might be close to a vee\footnote{If one blows up to a vee solution, an argument using Weiss' monotonicity formula as in \cite[Lemma 4.4]{CFFS25} yields that $u$ had to be a vee to begin with. Instead, to obtain quantitative estimates as in Corollary~\ref{cor:minimizers_regularity}, one needs to rule out nonflat global solutions blowing down to a vee profile.}. In order to obtain a result analogous to the boundedness of second fundamental form of free boundaries in Corollary~\ref{cor:minimizers_regularity}, one needs to take full advantage of the stability inequalities \eqref{eq:stab1}-\eqref{eq:stab2} in order to deduce that global stable solutions have flat free boundaries, even if they can blow-down to a vee. 

\subsection{Stability of cones} 

The picture, as in the case of minimizers, changes in dimension $7$: the De Silva--Jerison cone is a nonflat $1$-homogeneous minimizing cone in $\R^7$, and therefore also a stable cone. In particular, stable solutions (or limits of classically stable solutions, see \cite{ESV24, DJS22}) do not necessarily have smooth free boundaries from dimension $d \ge 7$.

For homogeneous cones, stability admits a particularly clean \emph{spectral} reformulation on the spherical link, and this point of view has become increasingly important. If $U$ is a $1$-homogeneous cone, the linearized equation around $U$ is \cite{CJK04, DJS22, JS15, fernandez2025linearized}
\[
\left\{
\begin{array}{rcll}
\Delta \varphi &=& 0 & \text{in}\quad \{U > 0\},\\
\varphi_\nu + H\varphi &=& 0 & \text{on}\quad \Reg(U),
\end{array}
\right.
\]
where $\Reg(U)$ denotes the regular part of $\partial\{U>0\}$, $\nu=\nabla U$ is the inward unit normal, and $H=-U_{\nu\nu}$. Following the tradition in minimal surface theory, we refer to the linearized equation as the Jacobi equation and its solutions as Jacobi fields. In order to study the behavior of Jacobi fields, by denoting
\[
\Omega_S:=\{U>0\}\cap \mathbb S^{d-1},\qquad \Gamma_S:=\Reg(U)\cap\mathbb S^{d-1},
\]
 one studies the Robin spectral problem on $\Omega_S$:
\[
\left\{
\begin{array}{rcll}
\Delta_{\mathbb S^{d-1}} \varphi &= &\lambda \varphi 
& \text{in }\Omega_S,\\
\partial_\nu \varphi + H\,\varphi& =& 0
& \text{on }\Gamma_S,
\end{array}
\right.
\]
where $\Delta_{\mathbb{S}^{d-1}}$ denotes the Laplace-Beltrami operator on $\mathbb{S}^{d-1}$.  With this sign convention, the principal eigenvalue $\lambda(U)$ is characterized variationally by
\[
-\lambda(U)
=
\inf_{\int_{\Omega_S}v^2>0} \frac{\displaystyle\int_{\Omega_S} |\nabla_{\mathbb S^{d-1}}v|^2 -\int_{\Gamma_S} H\,v^2} {\displaystyle\int_{\Omega_S} v^2},
\]
where the infimum is taken over $v\in C_c^\infty(\mathbb S^{d-1}\setminus\Sing(U))$. The stability condition is equivalent to the threshold inequality
\[
\lambda(U) \le \Big(\frac{d-2}{2}\Big)^2.
\]
Moreover, for a spherical eigenvalue $\lambda$, homogeneous Jacobi fields $r^\alpha\varphi(\theta)$ have degrees
\[
\alpha_\pm
=
-\frac{d-2}{2}\pm \sqrt{\Big(\frac{d-2}{2}\Big)^2-\lambda},
\]
so the spectrum directly controls the asymptotic rates of perturbations of the cone (see \cite{DJS22}). In the decay notation $r^{-\gamma}\varphi$ used in \cite{fernandez2025linearized}, one has $\gamma=-\alpha$. 

Recent work of Firester--Tsiamis--Wang pushes this spectral viewpoint in two directions. In \cite{firester2026stability} they prove strict stability inequalities for the cohomogeneity-one one-phase cones with bi-orthogonal symmetry (analogues of Lawson cones in minimal surfaces), showing that this family is strictly stable from dimension \(d = 7\) onward (instability for $d\le 6$ being shown in \cite{hong2015singular}); in particular, they obtain quantitative upper bounds for the principal Robin eigenvalue $\lambda(U)$ and corresponding bounds for the homogeneity degrees \(\alpha_\pm\), which yields decay estimates for Jacobi fields and applications to generic regularity (see \cite{fernandez2025linearized}). In the companion work \cite{firester2026area}, they develop the connection with capillary minimal hypersurfaces, construct nonflat capillary cones with bi-orthogonal symmetry for every dimension and contact angle, and use this correspondence to produce further examples of singular minimizing free boundaries for the Alt--Caffarelli functional. These are the first examples of nonsmooth minimizers for the Bernoulli problem outside the De Silva--Jerison cone in Theorem~\ref{thm:desilvajerison}. 

\subsection{A regularity conjecture} 
The previous observations lead to the stable analogue of the classical ``critical dimension $7$'' picture:

\begin{conj}
Every global classical stable solution to the one-phase problem in $\R^d$ has flat free boundary for $d\le 6$.
Equivalently, in dimensions $d\le 6$, classical stable solutions satisfy interior curvature bounds depending only on the dimension and the distance to the boundary.
\end{conj}

The conjecture is true in dimension $d=2$ (by a logarithmic cut-off argument). It is false in dimensions $d\ge 7$, because singular minimizing cones exist there, and one can construct families of classical minimizing solutions approximating them (as done in \cite{DJS22, ESV24}). Very recently it was proved to hold in dimension $d=3$, see \cite{CFFS25}; and in dimension $d =4$, see \cite{fernandezreal2026stability}. More precisely, in \cite{CFFS25, fernandezreal2026stability} it is shown that any global classical stable solution in $\R^d$ with $d \in \{3, 4\}$ is flat in the sense that
\[
D^2u\equiv 0 \quad\text{in }\{u>0\},
\]
so the free boundary consists of either one hyperplane or two parallel hyperplanes; as a consequence, one gets local curvature estimates for classical stable solutions in dimensions $3$ and $4$: 

\begin{thm} [\cite{CFFS25, fernandezreal2026stability}]
\label{thm:CFFS}
Let $u$ be a classical stable solution to the one-phase problem in $B_1\subset \R^d$, with $0\in \partial\{u > 0\}$ and $d\in \{3, 4\}$. Then, 
    \[
    |D^2 u |\le C \quad\text{in}\quad \{u > 0\}\cap B_{1/2}, 
    \]
    for some universal constant $C$. 
    In particular, the second fundamental form of the free boundary $|{\rm II}_{\Gamma_u}|$ is bounded by a universal constant in $B_{1/2}$.
\end{thm}

\section{Solutions beyond stability}

Once one leaves the minimizing or stable regimes, the one-phase problem becomes substantially more flexible. This reflects the nonconvexity of the Alt--Caffarelli functional: the class of critical points is much larger than the class of minimizers, and some of the key tools of the minimizing theory, such as nondegeneracy or a direct improvement-of-flatness mechanism, are no longer available in full generality. Thus, beyond minimizers one should think less in terms of ``regularity from energy'' and more in terms of the precise notion of solution being used, the topology of the positivity set, and the stability of the class under blow-up and compactness arguments \cite{Caf88,HHP11,CFFS25,kriventsov2025rectifiability, fernandezreal2026stability}.

The general theory splits into several complementary viewpoints. The first is the study of \emph{entire classical solutions}, where the free boundary is smooth, and the problem can be rephrased as an overdetermined problem on the positivity set. This leads to the language of \emph{exceptional domains} and connects the subject to minimal surfaces, Serrin-type problems, and fluid mechanics
\cite{HHP11,Tra14,khavinson2013overdetermined}. The second is the study of \emph{weak solutions}, where the Bernoulli condition is encoded through the distributional Laplacian and the free-boundary measure. This formulation is particularly close to (low regularity) Serrin-type overdetermined problems and is well-suited to local regularity questions. The third viewpoint is the study of \emph{closed classes of limits} of classical solutions, where one
sacrifices pointwise regularity in exchange for compactness and
measure-theoretic structure. Recent work of Kriventsov--Weiss shows that this latter direction is much more rigid than one might have expected a priori \cite{kriventsov2025rectifiability}.

These viewpoints overlap but should not be conflated. In the Alt--Caffarelli problem, minimizers are weak solutions  and weak solutions (in the Alt--Caffarelli definition) are stationary, but stationary solutions form a strictly larger class
\cite{Wei98,domingo2026regularity}. On the other hand, the weak class is not the maximal compactness class: limits of classical solutions may fail to remain weak, which is precisely one reason for introducing the variational/stationary framework of Kriventsov--Weiss \cite{kriventsov2025rectifiability}.

\subsection{Notions of solution}

\begin{figure}[t]
    \includegraphics[scale = 0.9]{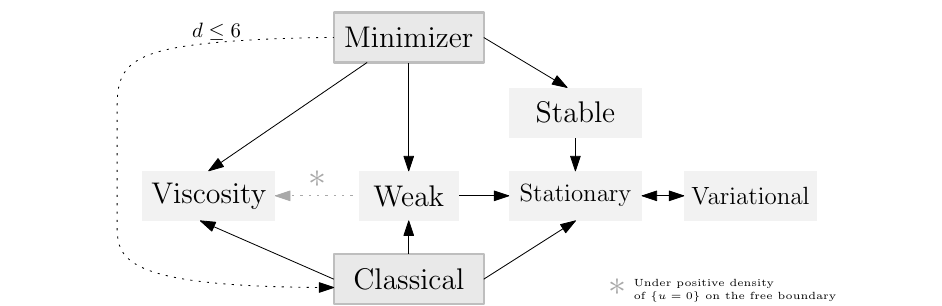}
    \centering
    \caption{\small Relationship between the different notions of solution. The notion of weak solution is to be understood in the Alt--Caffarelli sense; see Subsection~\ref{sssec:weak_sol}.}
    \label{fig:diagram}
\end{figure}

For the purposes of a global theory, it is important to distinguish several notions of solutions. We refer to Figure~\ref{fig:diagram} for a diagram representing the relationship between the notions we have introduced in this overview. 

\subsubsection{Classical solutions} At the most geometric level, one has \emph{classical solutions}: nonnegative functions
$u$ such that $\Omega_u:=\{u>0\}$ has smooth boundary\footnote{We say that a set has smooth boundary if the set can be locally expressed as the subgraph of a smooth function around each boundary point. In particular, the set cannot have boundary points with density one.} and
\[
\Delta u = 0 \quad \text{in } \Omega_u, \qquad u=0,\ |\nabla u|=1 \quad \text{on } \Gamma_u:=\partial \Omega_u.
\]
For such solutions, the positivity set itself is the main unknown, and one can bring ideas from differential geometry and overdetermined boundary value problems. While they are the easiest type of solution to work with from a regularity point of view, they present challenges: they do not form a closed class, and this is   due to the lack of density estimates for the zero set and to the lack of curvature estimates for the free boundary. They can, however, be used as the main elements of a relevant class, which can then be expanded by compactness (such as for stable or stationary solutions). 

\subsubsection{Viscosity solutions}
A second class is given by \emph{viscosity solutions}. This is the natural nonvariational
framework introduced by Caffarelli, designed to capture solutions without assuming
\emph{a priori} smoothness of the free boundary and to be stable under uniform limits
\cite{Caf88,CS05}. It is therefore well suited to compactness arguments (blow-up,
blow-down) and to the study of singular or rough interfaces.

A nonnegative continuous function $u$ is a \emph{viscosity solution} of the one-phase problem if
\[
\Delta u=0 \quad \text{in } \Omega_u:=\{u>0\}
\]
and the Bernoulli condition is encoded through comparison with smooth test functions
touching $u$ at free boundary points. More precisely, any strict comparison subsolution (resp. supersolution) cannot touch $u$ from below (resp. from above) at any free boundary point $x_0\in \Gamma_u$. Here, we say that $\phi$ is a (strict) comparison subsolution if $\phi\in C^2(\{\phi > 0\})$, $\Delta \phi \ge 0$ in $\{\phi > 0\}$, $\partial\{\phi > 0\}$ is a $C^2$ manifold  and $|\nabla \phi| > 1$ on $\partial\{\phi > 0\}$. (Resp., a comparison supersolution satisfies $\Delta \phi \le 0$ in $\{\phi > 0\}$ and $|\nabla \phi| < 1$ on $\partial\{\phi > 0\}$.)

This is the free-boundary analogue of the viscosity formulation for fully nonlinear PDEs (see \cite{CaffarelliCabre1995, fernandez2022regularity}) and provides a robust notion of solution even when $\Gamma_u$ has singularities
\cite{Caf88,CS05}.

The viscosity class is broad enough to include many non-minimizing configurations, and, crucially, it is closed under locally uniform convergence. This makes it the natural setting for limit problems and for local regularity arguments based on compactness and barrier
constructions. In regimes where one has \emph{flatness} (or, more generally, an appropriate nondegeneracy and density control), viscosity solutions enjoy a powerful local theory: Harnack inequalities, improvement of flatness, and $C^{1,\alpha}$ regularity of the free boundary at regular points \cite{Caf87,Caf88,CS05,de2011free}.

Precisely because it is nonvariational and very flexible, the viscosity class is also very broad. Properties that are automatic for minimizers (or for variational solutions) may fail for general
viscosity solutions unless they are imposed as hypotheses. For instance, nondegeneracy is not automatic: one can have solutions with very slow growth away from $\Gamma_u$, and the free boundary may contain components that are invisible at the level of a naive blow-up. Similarly, the viscosity class admits global solutions with complicated topology or large-scale geometry, so one should not expect a clean global classification without additional structure (e.g.\ monotonicity, graphicality, finite index, or quantitative flatness assumptions); see, for example, \cite{BSS76, EFeY23, basulto2024one}. 

  \subsubsection{Variational and stationary solutions}

A natural non-minimizing notion at the variational level is that of stationarity under inner variations: it corresponds to the criticality condition of the functional under perturbations of the form $u_t = u\circ ({\rm Id}+t X)$. This gives rise to the notion of \emph{stationary solutions} (as seen in \eqref{eq:stationary_def} in Section~\ref{sec:stable}). 

For compactness, however, it is useful to allow a relaxed phase. Following Weiss and
Kriventsov--Weiss \cite{Wei98,kriventsov2025rectifiability}, one considers pairs
$(u,\chi)$, where $u\ge0$ is Lipschitz in an open set $\Omega\subset \R^d$, harmonic in $\{u>0\}$, and $\chi:\Omega \to \{0, 1\}$ is Borel measurable and satisfies 
\[
\chi_{\{u>0\}}\leq \chi
\qquad\text{a.e. in }\Omega.
\]
The pair is called a variational solution if
\[
\int_\Omega \Bigl(   (|\nabla u|^2+\chi)\, {\rm div} X  -2\nabla u\cdot DX\nabla u \Bigr)=0
\qquad
\text{for every }X\in C_c^{0,1}(\Omega;\R^d).
\]
This is the first-variation identity for the relaxed energy
$\int_\Omega(|\nabla u|^2+\chi)$ under domain deformations. If one has the sharper identity with $\chi=\chi_{\{u>0\}}$, then $u$ is a stationary solution, as defined above.

The main result of Kriventsov--Weiss shows that the relaxation is only a compactness device: apart from the degenerate zero solution, no excess phase remains.

\begin{thm}[{\protect \cite{kriventsov2025rectifiability}}]\label{thm:kriventsov-weiss}
Let $(u,\chi)$ be a \emph{variational solution} of the Bernoulli problem in an open connected set $\Omega\subset\R^d$. Then either
$u\equiv0$ or
$
\chi=\chi_{\{u>0\}}
$ a.e. in $\Omega$.  Hence, every nontrivial variational solution is stationary. Moreover, $\partial\{u>0\}$ is countably $\cH^{d-1}$-rectifiable, has locally finite $\cH^{d-1}$ measure, and the measure $\Delta u$ can be identified completely. Limits of uniformly Lipschitz classical solutions and of suitable singular perturbation problems belong to this class.
\end{thm}

Thus, the relaxed formulation is closed under limits, while the theorem recovers the classical sharp phase a posteriori. The theorem therefore gives a robust framework for non-minimizing critical points: even without minimality or stability, the topological free boundary cannot be an arbitrary rough set, but is forced to carry codimension-one geometric measure-theoretic structure.

\subsubsection{Weak solutions}
\label{sssec:weak_sol}

There is also a more measure-theoretic weak formulation, going back to Alt--Caffarelli and refined by Weiss,  \cite{AC81,Wei98,domingo2026regularity}, in which the Bernoulli condition is encoded through the distributional Laplacian rather than imposed pointwise on a smooth free boundary. We say that $u$ is a \emph{weak solution} in $B\subset \R^d$ (in the sense of Alt--Caffarelli) if:
\begin{itemize}
  \item $u\in C(B)$, $u\ge0$, and $-\Delta u=0$ in the open set $\{u>0\}$;
  \item for every $D\Subset B$ there exist constants $0<c_D\le C_D<\infty$ such
  that
  \[
  c_D\le
  \frac1r\fint_{\partial B_r(x)}u\,d\cH^{d-1}   \le C_D
  \]
  whenever $x\in\partial\{u>0\}$ and $B_r(x)\subset D$;

  \item the identity
  \[
  -\int_B \nabla u\cdot \nabla \eta\,dx =
  \int_{\partial_{\mathrm{red}}\{u>0\}\cap B}  \eta\,d\cH^{d-1}
  \]
  holds for every $\eta\in C_c^\infty(B)$.
\end{itemize}
The second condition is a two-sided linear growth assumption: the lower bound is the nondegeneracy of $u$ at the free boundary, while the upper bound is  Lipschitz control. The third condition says that
the singular part of $\Delta u$ is precisely the free-boundary measure $ \cH^{d-1}\res\partial_{\rm red}\{u>0\}$, where $\partial_{\rm red}\{u>0\}$ denotes the reduced boundary of $\{u > 0\}$. The first two conditions imply that $\{u>0\}$ has locally finite perimeter, so its reduced boundary is well defined \cite[Remark 2.4]{domingo2026regularity}. Thus, the definition already builds in a rectifiable free-boundary measure; this is a feature of the Alt--Caffarelli weak theory, not of every distributional formulation of the Bernoulli condition.

There is, however, a closely related broader source of weak formulations. If $\Omega$ is a rough domain and $G$ is a Green function, extended by zero outside $\Omega$, then the boundary term in the distributional Laplacian is the harmonic measure. When this measure is absolutely continuous with respect to the surface measure, its density is the Poisson kernel. Thus, on smooth boundaries, the Bernoulli condition can be viewed as the prescription of the relevant Poisson kernel. This is one of the links between one-phase free boundaries, harmonic-measure regularity, and Serrin-type questions \cite{KenigToro99,FZ26,domingo2026regularity}. These harmonic-measure weak formulations are slightly broader than the Alt--Caffarelli notion of weak solution, and should not be automatically identified with stationary Alt--Caffarelli critical points; stationarity requires additional measure-theoretic structure.

From this perspective, (Alt--Caffarelli) weak solutions provide a convenient bridge between
Alt--Caffarelli free boundaries and Serrin-type overdetermined problems. Motivated by
the work of Figalli and Zhang on Serrin's problem in rough domains \cite{FZ26},
Domingo-Pasarin and Ros-Oton \cite{domingo2026regularity} prove that the weak formulation is strong enough to
recover smoothness once the free boundary is Lipschitz. In a way,  weak solutions are well adapted to local regularity, but they are not
the maximal compactness class, since the vee solution $|x\cdot e|$ is not a weak solution (even if it is the limit of classical solutions). In this sense, weak solutions occupy an
intermediate position between stationarity and minimality: they are stronger than mere stationarity, but flexible enough
to include the natural weak formulations arising from overdetermined boundary value
problems.

\subsection{Exceptional domains and a zoo of examples}

For a global classical solution, the positivity set $\{u > 0\}$ is called an \emph{exceptional
domain}: it carries a positive harmonic function $u$ vanishing on $\Gamma_u$
and satisfying the overdetermined condition $|\nabla u|=1$ there; see \cite{HHP11,khavinson2013overdetermined,Tra14}.

As we have seen, the half-space solution
\[
u(x) = (x\cdot e)_+, \qquad e \in \mathbb{S}^{d-1},
\]
is the canonical example, which also includes the family of solutions 
\[
u(x) = (x\cdot e-a)_+ + (x\cdot e - b)_-\qquad\text{for}\quad e\in \mathbb{S}^{d-1}
\]
for $a > b$ converging to a vee as $b \uparrow a$. In particular, complements of strips  of the form 
\[
\R^d\setminus \{x : b \le x\cdot e \le a\}\qquad\text{for some}\quad e\in \mathbb{S}^{d-1}, -\infty < b < a \le \infty,
\]
are exceptional domains. 

Similarly, after the appropriate normalization, the exterior of a ball is an exceptional
domain. In the plane this gives the disk-complement solution
\[
u(x) = (\log |x|)_+,
\]
and in higher dimensions $d \ge 3$ one has the corresponding asymptotic behavior according to the fundamental solutions of the Laplacian: 
\[
u(x) = \frac{1}{d-2}(1- |x|^{2-d})_+.
\]
In fact, this is the only case in which an entire solution has a bounded free boundary (see, e.g., \cite{reichel1997radial, FFS26}). Thus, even among the most elementary examples, one already sees free boundaries with
different topology and different large-scale behavior.

Beyond these basic configurations, Hauswirth--H\'elein--Pacard constructed the
double-hairpin solution in the plane, thus discovering a catenoid-type
family of exceptional domains \cite{HHP11}. The zero-level set (after a rescaling, translation, and rotation) is given by 
\[
\left\{(x_1, x_2)\in \R^2 : |x_2|\ge \tfrac{\pi}{2}+ \cosh(x_1)\right\}.
\]
In higher dimensions, Liu--Wang--Wei
constructed smooth higher-dimensional catenoid-type solutions (satisfying the assumptions in \cite{FFS26}, and thus with expansions like \eqref{eq:gasympt}), showing that the
phenomenon is by no means restricted to $d=2$ \cite{liu2021smooth}. 

There are also periodic
families of planar solutions, originating in the fluid mechanical work of
Baker--Saffman--Sheffield 
\cite{BSS76} (see also \cite{crowdy2011analytical}). These examples are conceptually important: they show that,
outside the minimizing regime, one should not expect the positivity set to resemble a
single half-space, even asymptotically.

In fact, the periodic examples already indicate how subtle the closure of the classical
class can be. As emphasized by Kriventsov--Weiss, suitable blow-downs of such planar
periodic solutions may converge to $\alpha |x_2|$ for any chosen $\alpha \in (0,1)$, showing
that even limits of classical solutions may fall outside the classical Bernoulli picture in
a pointwise sense \cite{kriventsov2025rectifiability}. By contrast, if the classical solution is assumed to be stable, then blow-downs of the form $\alpha |x_2|$ can only occur with $\alpha=1$ \cite{CFFS25}. This illustrates why the choice of solution class is so important.

The double-hairpin solution also leads to periodic examples in higher dimensions. Basulto--Kamburov \cite{basulto2026periodic} construct a family of global classical solutions in $\R^3$ bifurcating from the product of the planar double hairpin with a line. These solutions depend nontrivially and periodically on the additional variable, and their free boundaries remain noncompact modulo this periodicity.

\subsubsection{Rigidity under graphicality}

Even if a general classification is out of reach, there are some situations in which a
strong global geometric assumption forces a non-minimizing solution back into the
minimizing world.

One example of such a result is the recent work \cite{EFeY23}. There, viscosity solutions are considered whose contact set is a subgraph of a continuous function, or more generally satisfies a monotonicity-type graphical condition. The main
result is that such solutions are automatically minimizers of the Alt--Caffarelli functional. Consequently, in dimensions $d\le d_\ast=7$, an entire solution whose contact set is a subgraph of a continuous function must be a half-space solution. This is a natural one-phase analogue of Bernstein-type rigidity for minimal graphs.

A precursor in the minimizing setting is the gradient estimate of De Silva and
Jerison \cite{de2011gradient}. They proved a one-phase analogue of the
Bombieri--De Giorgi--Miranda interior gradient bound for minimal graphs:
roughly speaking, if a solution is monotone in one direction and its free
boundary is graphical, then one obtains an interior Lipschitz estimate for the
graph under a quantitative geometric assumption on the positive phase. In
particular, this applies to energy-minimizing solutions and shows that every
energy-minimizing graphical free boundary is smooth.  

\subsection{Homogeneous solutions}

As in the minimizing and stable theories, homogeneous solutions play a central role,
since they arise naturally in blow-up and blow-down procedures. However, once one is no
longer restricted to minimizing or stable variational classes, the homogeneous world is much richer and much less
understood.

A one-homogeneous solution in $\R^d$ corresponds to an overdetermined eigenvalue
problem on a domain in $\mathbb{S}^{d-1}$, so the classification of homogeneous solutions can
be reformulated as a shape-optimization problem on the sphere. In low dimensions,
classical examples were for a long time extremely scarce. In dimension $3$, besides the
half-space solution, the main known example was the rotationally symmetric
Alt--Caffarelli solution (the three-dimensional version of the De Silva--Jerison cone). Recent work of Hines--Kolesar--McGrath \cite{hines2025new} changes this picture: they construct, for all sufficiently large integers $k$, new homogeneous
solutions in $\R^3$ whose spherical links have $k$ boundary components, and they also
construct an infinite family of homogeneous solutions in dimension $4$. In particular, even the low-dimensional homogeneous theory is considerably richer.

By contrast, Parise--Zhu prove that a $1$-homogeneous Bernoulli solution in $\R^3$ that is classical away from the origin and whose spherical positivity set is a smooth annulus is, up to rotation, the axially symmetric Alt--Caffarelli cone \cite[Theorem~1.6]{parise2026uniqueness}.

In higher dimensions, Firester--Tsiamis--Zhao construct homogeneous one-phase solutions from isoparametric foliations of the sphere, obtaining new examples in every dimension $d\ge5$ and rigidity within the corresponding isoparametric classes \cite{firester2026isoparametric}. Stability and minimality are separate questions for these examples.

This is an important warning for the general classification problem, since homogeneous
solutions are the first objects one encounters in blow-up analysis, so any attempt at a
global rigidity theorem must eventually confront the question of which homogeneous
limits are possible.

\subsection{Finite Morse index}

A natural intermediate class between stable solutions and general classical solutions is given by critical points with finite index. If $u$ is a classical solution in a domain $D$, with positivity set $\Omega_u:=\{u>0\}$ and free boundary $\Gamma_u:=\partial\{u>0\}$, its second variation is encoded by the quadratic form 
\[
\mathcal Q_D(\xi,\xi)
= \int_{D\cap\Omega_u} |\nabla \xi|^2\,dx - \int_{D\cap\Gamma_u} H\,\xi^2\,d\cH^{d-1},
\]
where $\nu=\nabla u$ is the inward unit normal to $\Omega_u$ and
$H=-u_{\nu\nu}$ denotes the mean curvature of $\Gamma_u$ with this convention. The Morse index of $u$ in $D$ is the maximal dimension of a subspace
on which $\mathcal Q_D$ is negative definite. Thus stability corresponds to index
zero, while finite Morse index means that the solution is unstable only in finitely
many independent directions.

This is the free-boundary analogue of the finite-index class for minimal surfaces. The
main heuristic is that finite index should imply \emph{asymptotic stability}: all
instability is confined to a compact core, while the ends of the solution recover the
rigidity of stable solutions. Basulto--Kamburov made this precise in the one-phase
setting. They prove that an entire classical solution of finite Morse index is stable
outside a compact set, and they identify the index with the number of negative
eigenvalues, counted with multiplicity, of the Jacobi--Robin problem
\[
-\Delta \phi=\mu\phi
\quad\text{in }\Omega_u, \qquad \partial_\nu\phi+H\phi=0 \quad\text{on }\Gamma_u,
\]
under exhaustion by bounded domains \cite{basulto2024one}.

In the plane this leads to a satisfactory classification. Basulto--Kamburov prove
that if a global classical solution in $\R^2$ has connected positivity set and finite
Morse index, then its free boundary has finite topology. Combining this with
Traizet's classification of planar exceptional domains, one obtains that such a
solution is, up to similarity, either the half-plane solution, the disk-complement
solution, or the double-hairpin solution. The first is stable, while the latter two
have Morse index one \cite{Tra14,basulto2024one}. This is the closest one-phase
analogue of the classical finite-index theory for complete minimal surfaces in
Euclidean space.

The very recent work \cite{FFS26} extends this picture to three dimensions.
More generally, it yields an analogous classification in any dimension
under the hypothesis that all classical stable solutions and all
homogeneous stable solutions are flat. Thanks to
\cite{fernandezreal2026stability}, the classification therefore also
holds in four dimensions.

In particular, solutions with finite Morse index in $\R^3$ are either
one-dimensional, have contact set equal to a ball, or are axially
symmetric. In the last case, after rescaling, translating, and rotating,
we may write
\[
u(x_1,x_2, x_3)=u(r,x_3),
\qquad
r=\sqrt{x_1^2+x_2^2}.
\]
The solution is even in $x_3$ and satisfies
\[
\partial_{x_3}u>0 \quad\text{in }\{x_3>0\}\cap\{u>0\}, \qquad \partial_r u<0 \quad\text{in }\{r>0\}\cap\{u>0\}.
\]
In this case the zero phase is a bounded-width catenoidal region:
\[
\{u=0\}
=
\{(x',x_3): |x'|\ge 1,\ |x_3|\le g(|x'|)\},
\]
where $g:[1,\infty)\to[0,\infty)$ satisfies $g(1)=0$, $g'(r)>0$ for $r>1$, and
$\lim_{r\to\infty}g(r)<\infty$. More precisely, for some $b,c>0$ and every
$\alpha\in(0,1)$,
\begin{equation}
\label{eq:gasympt}
g(r)=b-c\,r^{-1}+O(r^{-2-\alpha}),
\end{equation}
with analogous expansions for higher derivatives; \cite{FFS26}. In fact, these asymptotics are shown more generally for flat solutions at infinity for the Bernoulli problem (or solutions in exterior domains with flat ends).

\addtocontents{toc}{\protect\setcounter{tocdepth}{-1}}

\section*{Acknowledgements}

The author is grateful to Max Engelstein, Dennis Kriventsov, and Hui Yu for their comments and suggestions on an earlier version of this manuscript.

\section*{AI declaration}

AI tools were used to correct spelling, improve exposition, and check the bibliography.

\section*{Conflict of interest declaration}

The author has no conflicts of interest to declare. 

\addtocontents{toc}{\protect\setcounter{tocdepth}{2}}

\begingroup
\renewcommand{\bibliofont}{\small}
\let\reviewbibitem\bibitem
\renewcommand{\bibitem}[2][]{\def\reviewkey{openai2026critical}\def\currentkey{#2}\ifx\currentkey\reviewkey\reviewbibitem[OAI26]{#2}\else\reviewbibitem[#1]{#2}\fi}
\bibliographystyle{alpha}
\bibliography{mybibliography}

@article{AC81,
  author  = {Alt, H. W. and Caffarelli, L. A.},
  title   = {Existence and regularity for a minimum problem with free boundary},
  journal = {J. Reine Angew. Math.},
  volume  = {325},
  year    = {1981},
  pages   = {105--144}
}

@article{AP86,
  author  = {Alt, H. W. and Phillips, D.},
  title   = {A free boundary problem for semilinear elliptic equations},
  journal = {J. Reine Angew. Math.},
  volume  = {368},
  year    = {1986},
  pages   = {63--107}
}

@article{audrito2022interface,
  author  = {Audrito, A. and Serra, J.},
  title   = {Interface regularity for semilinear one-phase problems},
  journal = {Adv. Math.},
  volume  = {403},
  year    = {2022},
  pages   = {108380}
}

@article{basulto2024one,
  author  = {Basulto, J. and Kamburov, N.},
  title   = {One-phase free boundary solutions of finite {M}orse index},
  journal = {J. Differential Equations},
  volume  = {410},
  year    = {2024},
  pages   = {319--345}
}

@book{BL82,
  author    = {Buckmaster, J. D. and Ludford, G. S. S.},
  title     = {Theory of Laminar Flames},
  series    = {Cambridge Monographs on Mechanics},
  publisher = {Cambridge Univ. Press},
  address   = {Cambridge},
  year      = {1982}
}

@article{BSS76,
  author  = {Baker, G. R. and Saffman, P. G. and Sheffield, J. S.},
  title   = {Structure of a linear array of hollow vortices of finite cross-section},
  journal = {J. Fluid Mech.},
  volume  = {74},
  year    = {1976},
  pages   = {469--476}
}

@article{Caf87,
  author  = {Caffarelli, L. A.},
  title   = {A {H}arnack inequality approach to the regularity of free boundaries. {I}. {L}ipschitz free boundaries are {$C^{1,\alpha}$}},
  journal = {Rev. Mat. Iberoam.},
  volume  = {3},
  year    = {1987},
  pages   = {139--162}
}

@article{Caf88,
  author  = {Caffarelli, L. A.},
  title   = {A {H}arnack inequality approach to the regularity of free boundaries. {III}. Existence theory, compactness, and dependence on {$X$}},
  journal = {Ann. Scuola Norm. Sup. Pisa Cl. Sci. (4)},
  volume  = {15},
  year    = {1988},
  pages   = {583--602}
}

@article{Caf89,
  author  = {Caffarelli, L. A.},
  title   = {A {H}arnack inequality approach to the regularity of free boundaries. {II}. {F}lat free boundaries are {L}ipschitz},
  journal = {Comm. Pure Appl. Math.},
  volume  = {42},
  year    = {1989},
  pages   = {55--78}
}

@book{CaffarelliCabre1995,
  author    = {Caffarelli, L. A. and Cabr{\'e}, X.},
  title     = {Fully Nonlinear Elliptic Equations},
  series    = {Amer. Math. Soc. Colloq. Publ.},
  volume    = {43},
  publisher = {Amer. Math. Soc.},
  address   = {Providence, RI},
  year      = {1995}
}

@article{CFFS25,
  author = {Chan, H. and Fern{\'a}ndez-Real, X. and Figalli, A. and Serra, J.},
  title  = {Global stable solutions to the free boundary {A}llen--{C}ahn and {B}ernoulli problems in 3{D} are one-dimensional},
  journal = {J. Amer. Math. Soc.},
  year   = {2026},
  pages     = {1027--1151}, 
  volume    = {39},
  number    = {4}
}

@incollection{CJK04,
  author    = {Caffarelli, L. A. and Jerison, D. and Kenig, C. E.},
  title     = {Global energy minimizers for free boundary problems and full regularity in three dimensions},
  booktitle = {Noncompact Problems at the Intersection of Geometry, Analysis, and Topology},
  series    = {Contemp. Math.},
  volume    = {350},
  publisher = {Amer. Math. Soc.},
  address   = {Providence, RI},
  year      = {2004},
  pages     = {83--97}
}

@book{CS05,
  author    = {Caffarelli, L. A. and Salsa, S.},
  title     = {A Geometric Approach to Free Boundary Problems},
  series    = {Grad. Stud. Math.},
  volume    = {68},
  publisher = {Amer. Math. Soc.},
  address   = {Providence, RI},
  year      = {2005}
}

@article{crowdy2011analytical,
  author  = {Crowdy, D. G. and Green, C. C.},
  title   = {Analytical solutions for von {K}{\'a}rm{\'a}n streets of hollow vortices},
  journal = {Phys. Fluids},
  volume  = {23},
  number  = {12},
  year    = {2011},
  pages   = {126602}
}

@article{CV95,
  author  = {Caffarelli, L. A. and V{\'a}zquez, J. L.},
  title   = {A free-boundary problem for the heat equation arising in flame propagation},
  journal = {Trans. Amer. Math. Soc.},
  volume  = {347},
  year    = {1995},
  pages   = {411--441}
}

@article{de2011free,
  author  = {{De Silva}, D.},
  title   = {Free boundary regularity for a problem with right hand side},
  journal = {Interfaces Free Bound.},
  volume  = {13},
  number  = {2},
  year    = {2011},
  pages   = {223--238}
}

@article{de2011gradient,
  author  = {{De Silva}, D. and Jerison, D.},
  title   = {A gradient bound for free boundary graphs},
  journal = {Comm. Pure Appl. Math.},
  volume  = {64},
  number  = {4},
  year    = {2011},
  pages   = {538--555}
}

@article{DJ09,
  author  = {{De Silva}, D. and Jerison, D.},
  title   = {A singular energy minimizing free boundary},
  journal = {J. Reine Angew. Math.},
  volume  = {635},
  year    = {2009},
  pages   = {1--21},
  doi     = {10.1515/CRELLE.2009.074}
}

@article{DJS22,
  author  = {{De Silva}, D. and Jerison, D. and Shahgholian, H.},
  title   = {Inhomogeneous global minimizers to the one-phase free boundary problem},
  journal = {Comm. Partial Differential Equations},
  volume  = {47},
  year    = {2022},
  pages   = {1193--1216}
}

@misc{domingo2026regularity,
  author = {Domingo-Pasarin, J. and Ros-Oton, X.},
  title  = {Regularity of {L}ipschitz free boundaries for weak solutions of {A}lt--{C}affarelli type problems},
  year   = {2026},
  note   = {Preprint, arXiv:2601.20493}
}

@article{DSV21,
  author  = {{De Philippis}, G. and Spolaor, L. and Velichkov, B.},
  title   = {Regularity of the free boundary for the two-phase {B}ernoulli problem},
  journal = {Invent. Math.},
  volume  = {225},
  year    = {2021},
  pages   = {347--394}
}

@article{EE19,
  author  = {Edelen, N. and Engelstein, M.},
  title   = {Quantitative stratification for some free-boundary problems},
  journal = {Trans. Amer. Math. Soc.},
  volume  = {371},
  year    = {2019},
  pages   = {2043--2072}
}

@article{EFeY23,
  author  = {Engelstein, M. and Fern{\'a}ndez-Real, X. and Yu, H.},
  title   = {Graphical solutions to one-phase free boundary problems},
  journal = {J. Reine Angew. Math.},
  volume  = {804},
  year    = {2023},
  pages   = {155--195}
}

@article{ESV20,
  author  = {Engelstein, M. and Spolaor, L. and Velichkov, B.},
  title   = {Uniqueness of the blowup at isolated singularities for the {A}lt--{C}affarelli functional},
  journal = {Duke Math. J.},
  volume  = {169},
  year    = {2020},
  pages   = {1541--1601}
}

@article{ESV24,
  author  = {Edelen, N. and Spolaor, L. and Velichkov, B.},
  title   = {A strong maximum principle for minimizers of the one-phase {B}ernoulli problem},
  journal = {Indiana Univ. Math. J.},
  volume  = {73},
  year    = {2024},
  pages   = {1061--1096}
}

@misc{fernandez2025linearized,
  author = {Fern{\'a}ndez-Real, X. and Yu, H.},
  title  = {Linearized equation and generic regularity in the {A}lt--{C}affarelli problem},
  year   = {2025},
  note   = {Preprint, arXiv:2510.18330}
}

@book{fernandez2022regularity,
  author    = {Fern{\'a}ndez-Real, X. and Ros-Oton, X.},
  title     = {Regularity Theory for Elliptic {PDE}},
  series    = {Zurich Lectures in Advanced Mathematics},
  volume    = {28},
  publisher = {EMS Press},
  address   = {Berlin},
  year      = {2022},
  doi       = {10.4171/ZLAM/28}
}

@misc{FFS26,
  author = {Fern{\'a}ndez-Real, X. and Florit-Simon, E. and Serra, J.},
  title  = {Finite index solutions to the {B}ernoulli problem in three dimensions are axially symmetric},
  year   = {2026},
  note   = {Preprint, arXiv:2605.07913}
}

@misc{firester2026area,
  author = {Firester, B. and Tsiamis, R. and Wang, Y.},
  title  = {Area-minimizing capillary cones},
  year   = {2026},
  note   = {Preprint, arXiv:2601.18794}
}

@misc{firester2026stability,
  author = {Firester, B. and Tsiamis, R. and Wang, Y.},
  title  = {Stability inequalities for one-phase cones},
  year   = {2026},
  note   = {Preprint, arXiv:2601.16966}
}

@article{fr2019dcds,
  author  = {Fern{\'a}ndez-Real, X. and Ros-Oton, X.},
  title   = {On global solutions to semilinear elliptic equations related to the one-phase free boundary problem},
  journal = {Discrete Contin. Dyn. Syst.},
  volume  = {39},
  number  = {12},
  year    = {2019},
  pages   = {6945--6959}
}

@article{FY23,
  author = {Fern{\'a}ndez-Real, X. and Yu, H.},
  title  = {Generic properties in free boundary problems},
  journal = {Amer. J. Math.},
  year   = {2023},
  note   = {To appear; arXiv:2308.13209}
}

@article{FZ26,
  author  = {Figalli, A. and Zhang, Y. R.-Y.},
  title   = {Serrin's overdetermined problem in rough domains},
  journal = {J. Eur. Math. Soc. (JEMS)},
  year    = {2025},
  note    = {Published online first},
  doi     = {10.4171/JEMS/1726}
}

@article{HHP11,
  author  = {Hauswirth, L. and H{\'e}lein, F. and Pacard, F.},
  title   = {On an overdetermined elliptic problem},
  journal = {Pacific J. Math.},
  volume  = {250},
  year    = {2011},
  pages   = {319--334}
}

@misc{hines2025new,
  author = {Hines, C. and Kolesar, J. and McGrath, P.},
  title  = {New homogeneous solutions for the one-phase free boundary problem},
  year   = {2025},
  note   = {Preprint, arXiv:2509.09409}
}

@article{hong2015singular,
  author  = {Hong, G.},
  title   = {The singular homogeneous solutions to one phase free boundary problem},
  journal = {Proc. Amer. Math. Soc.},
  volume  = {143},
  number  = {9},
  year    = {2015},
  pages   = {4009--4015}
}

@article{JS15,
  author  = {Jerison, D. and Savin, O.},
  title   = {Some remarks on stability of cones for the one-phase free boundary problem},
  journal = {Geom. Funct. Anal.},
  volume  = {25},
  year    = {2015},
  pages   = {1240--1257}
}

@article{kamburov2022nondegeneracy,
  author  = {Kamburov, N. and Wang, K.},
  title   = {Nondegeneracy for stable solutions to the one-phase free boundary problem},
  journal = {Math. Ann.},
  volume  = {388},
  year    = {2024},
  pages   = {2705--2726}
}

@article{khavinson2013overdetermined,
  author  = {Khavinson, D. and Lundberg, E. and Teodorescu, R.},
  title   = {An overdetermined problem in potential theory},
  journal = {Pacific J. Math.},
  volume  = {265},
  number  = {1},
  year    = {2013},
  pages   = {85--111}
}

@article{KL18,
  author  = {Kriventsov, D. and Lin, F.},
  title   = {Regularity for shape optimizers: the nondegenerate case},
  journal = {Comm. Pure Appl. Math.},
  volume  = {71},
  year    = {2018},
  pages   = {1535--1596}
}

@article{KL19,
  author  = {Kriventsov, D. and Lin, F.},
  title   = {Regularity for shape optimizers: the degenerate case},
  journal = {Comm. Pure Appl. Math.},
  volume  = {72},
  year    = {2019},
  pages   = {1678--1721}
}

@misc{kriventsov2025min,
  author = {Kriventsov, D. and Weiss, G. S.},
  title  = {A min-max variational approach to the existence of gravity water waves},
  year   = {2025},
  note   = {Preprint; arXiv:2508.12159}
}

@article{kriventsov2025rectifiability,
  author  = {Kriventsov, D. and Weiss, G. S.},
  title   = {Rectifiability, finite {H}ausdorff measure, and compactness for non-minimizing {B}ernoulli free boundaries},
  journal = {Comm. Pure Appl. Math.},
  volume  = {78},
  number  = {3},
  year    = {2025},
  pages   = {545--591}
}

@article{liu2021smooth,
  author  = {Liu, Y. and Wang, K. and Wei, J.},
  title   = {On smooth solutions to one phase-free boundary problem in {$\mathbb{R}^n$}},
  journal = {Int. Math. Res. Not. IMRN},
  volume  = {2021},
  number  = {20},
  year    = {2021},
  pages   = {15682--15732}
}

@article{PY07,
  author  = {Petrosyan, A. and Yip, N. K.},
  title   = {Nonuniqueness in a free boundary problem from combustion},
  journal = {J. Geom. Anal.},
  volume  = {18},
  year    = {2008},
  pages   = {1098--1126}
}

@incollection{sirakov2002overdetermined,
  author    = {Sirakov, B.},
  title     = {Overdetermined elliptic problems in physics},
  editor    = {Berestycki, H. and Pomeau, Y.},
  booktitle = {Nonlinear {PDE}'s in Condensed Matter and Reactive Flows},
  series    = {NATO Sci. Ser. C Math. Phys. Sci.},
  volume    = {569},
  publisher = {Springer},
  address   = {Dordrecht},
  year      = {2002},
  pages     = {273--295},
  doi       = {10.1007/978-94-010-0307-0_11}
}

@article{Tra14,
  author  = {Traizet, M.},
  title   = {Classification of the solutions to an overdetermined elliptic problem in the plane},
  journal = {Geom. Funct. Anal.},
  volume  = {24},
  year    = {2014},
  pages   = {690--720}
}

@book{Vel23,
  author    = {Velichkov, B.},
  title     = {Regularity of the One-Phase Free Boundaries},
  series    = {Lecture Notes Unione Mat. Ital.},
  volume    = {28},
  publisher = {Springer},
  address   = {Cham},
  year      = {2023}
}

@article{Wei98,
  author  = {Weiss, G. S.},
  title   = {Partial regularity for weak solutions of an elliptic free boundary problem},
  journal = {Comm. Partial Differential Equations},
  volume  = {23},
  year    = {1998},
  pages   = {439--455}
}

@article{Wei99,
  author  = {Weiss, G. S.},
  title   = {Partial regularity for a minimum problem with free boundary},
  journal = {J. Geom. Anal.},
  volume  = {9},
  year    = {1999},
  pages   = {317--326}
}

@article{KenigToro99,
  author  = {Kenig, C. E. and Toro, T.},
  title   = {Free boundary regularity for harmonic measures and {P}oisson kernels},
  journal = {Ann. of Math. (2)},
  volume  = {150},
  number  = {2},
  year    = {1999},
  pages   = {369--454}
}

@misc{engelstein2025asymptotic,
  author = {M. Engelstein and D. Restrepo and Z. Zhao},
  title  = {On the asymptotic properties of solutions to one-phase free boundary problems},
  year   = {2025},
  note   = {Preprint, arXiv:2511.08393}
}

@article{edelen2024symmetric,
  author  = {N. Edelen and L. Spolaor and B. Velichkov},
  title   = {The symmetric (log-)epiperimetric inequality and a decay-growth estimate},
  journal = {Calc. Var. Partial Differential Equations},
  volume  = {63},
  number  = {1},
  year    = {2024},
  pages   = {Paper No. 2, 29}
}

@article{phillips1983minimization,
  title={A minimization problem and the regularity of solutions in the presence of a free boundary},
  author={Phillips, D.},
  journal={Indiana Univ. Math. J.},
  volume={32},
  number={1},
  pages={1--17},
  year={1983},
  doi={10.1512/iumj.1983.32.32001}
}

@article{reichel1997radial,
  title={Radial symmetry for elliptic boundary-value problems on exterior domains},
  author={Reichel, W.},
  journal={Arch. Ration. Mech. Anal.},
  volume={137},
  number={4},
  pages={381--394},
  year={1997},
  publisher={Springer}
}

@misc{fernandezreal2026stability,
  author = {Fern{\'a}ndez-Real, X. and Serra, J.},
  title  = {The stability conjecture for the {B}ernoulli problem in dimension four},
  year   = {2026},
  note   = {Preprint, arXiv:2609.05669}
}

@misc{carducci2026minimizing,
  author = {Carducci, M. and Velichkov, B.},
  title  = {A minimizing cone from {Cartan} polynomials in the {Alt-Caffarelli} problem},
  year   = {2026},
  note   = {Preprint, arXiv:2609.35179}
}

@misc{dephilippis2026optimal,
  author = {{De Philippis}, G. and Hirsch, J. and Nahon, M.},
  title  = {Optimal regularity for vectorial, higher order and non-minimizing {Bernoulli} problems},
  year   = {2026},
  note   = {Preprint, arXiv:2609.12180}
}

@misc{parise2026uniqueness,
  author = {Parise, D. and Zhu, J. J.},
  title  = {Uniqueness of free boundary minimal annuli},
  year   = {2026},
  note   = {Preprint, arXiv:2609.17512}
}

@misc{firester2026isoparametric,
  author = {Firester, B. and Tsiamis, R. and Zhao, Z.},
  title  = {Singularity models for the {Bernoulli} free boundary problem from isoparametric hypersurfaces},
  year   = {2026},
  note   = {Preprint, arXiv:2609.31617}
}

@misc{chan2026stableallen,
  author = {Chan, H. and Fern{\'a}ndez-Real, X. and Figalli, A. and Florit-Simon, E. and Serra, J.},
  title  = {Stable solutions of the {Allen--Cahn} equation in dimension three are one-dimensional},
  year   = {2026},
  note   = {Preprint, arXiv:2609.30194}
}

@misc{liu2026stabledegiorgi,
  author = {Liu, Y. and Luo, T. and Wang, K. and Wei, J. and Wei, Y. and Wu, K.},
  title  = {Stable {De Giorgi} conjecture of the {Allen--Cahn} equation in $\mathbb{R}^3$},
  year   = {2026},
  note   = {Preprint, arXiv:2609.21680}
}

@misc{basulto2026periodic,
  author = {Basulto, J. and Kamburov, N.},
  title  = {Nontrivial periodic solutions to the one-phase free boundary problem},
  year   = {2026},
  note   = {Preprint, arXiv:2610.05257}
}

@misc{openai2026critical,
  author = {{OpenAI}},
  title  = {The critical dimension for one-phase {Bernoulli} minimizers},
  year   = {2026},
  howpublished = {OpenAI Math Release preprint},
  note   = {\begingroup\def\UrlBreaks{\do\/\do\-}\url{https://github.com/openai/math/blob/main/preprints/The-critical-dimension-for-one-phase-Bernoulli-minimizers-September-24-2026/The-critical-dimension-for-one-phase-Bernoulli-minimizers-September-24-2026.pdf}\endgroup}
}
\endgroup

\end{document}